\documentclass[12pt,twoside]{article}
\usepackage[all]{xy}
\usepackage{a4,amsmath,amssymb,amsfonts,amscd,mathrsfs}
\reversemarginpar
\newcounter{Abschnitt}[section]

\newtheorem{theorem}{Theorem}[section]
\newtheorem{lemma}[theorem]{Lemma}
\newtheorem{corollary}[theorem]{Corollary}
\newtheorem{proposition}[theorem]{Proposition}

\newcommand{\Ext}{{\rm Ext}}

\newcommand{\Hom}{{\rm Hom}}

\newcommand{\ol}[1]{\overline{#1}}
\newcommand{\Pic}{{\rm Pic}}

\newcommand{\Quot}{{\rm Quot}}

\newcommand{\SU}{{\rm SU}}

\newcommand{\Sym}{{\rm Sym}}

\newcommand{\U}{{\rm U}}

\newcommand{\coker}{{\rm coker}}

\newcommand{\ext}{{\rm ext}}

\newcommand{\rk}{{\rm rk}}

\newcommand{\Ocal}{{\cal O}}

\newcommand{\ndop}{{\mathbb N}}

\newcommand{\pdopi}{{{\mathbb P}^1}}
\newcommand{\pdop}{{\mathbb P}}
\newcommand{\qdop}{{\mathbb Q}}

\newcommand{\zdop}{{\mathbb Z}}

\newcommand{\dual}{^\lor}

\newcommand{\inv}{^{-1}}

\newcommand{\rarpa}[1]{\stackrel{#1}{\rightarrow}}

\newcommand{\qed}{{ \hfill $\Box$}}

\title{Raynaud vector bundles}
\author{Georg Hein}
\date{June 27, 2007}

\begin{document}
\maketitle

\begin{abstract}
We construct vector bundles $R^\rk_\mu$ on a smooth projective curve $X$
having the property that for all sheaves $E$ of slope $\mu$ and rank
$\rk$ on $X$ we have an equivalence: $E$ is a semistable vector bundle $\iff$
$\Hom(R^\rk_\mu,E)=0$.

As a byproduct of our construction we obtain
effective bounds on $r$ such that the linear system $|R \cdot \Theta|$
has base points on $\U_X(r,r(g-1))$.
\end{abstract}

\section{Introduction}
{\bf Notation:}
Throughout this paper $X$ is a smooth projective curve over of genus $g$
over some algebraically closed field $k$.

Raynaud constructed in his article \cite{Ray} vector bundles $\{P_m\}_{m
\geq 1}$ with the property that $\mu(P_m)=\frac{g}{m}$ and $h^0(P_m
\otimes L) \ne 0$ for all linebundles $L$ of degree zero. We showed in
our article \cite{Hei} that the converse also holds, that is:
$h^0(E \otimes L) \ne 0$ for all linebundles of degree zero if and only
if, we have morphisms $P_{\rk(E)g+1} \to E$. Furthermore, Raynaud showed
that a vector bundle $E$ of rank two and slope $g-1$ is semistable if and
only if, there exists a linebundle $L$ of degree zero with
$h^0(E \otimes L)=0$. Thus, we deduced:
\begin{theorem}\label{intro-1}
For a coherent sheaf $E$ on $X$ of rank two and slope $g-1$ we have an
equivalence
\[ E \mbox{ is a semistable vector bundle } \iff \Hom(P_{2g+1},E)=0 \,. \]
\end{theorem}
This way we obtain another equivalent condition to semistability.
This condition is be very convenient, because we have to check the
behavior of $E$ only with respect to one bundle to decide, whether it is
semistable or not.
This motivates the following definition:\\

\begin{definition}
A vector bundle $R^\rk_\mu$ is called a Raynaud bundle, if we have an
equivalence 
\[ E \mbox{ is a semistable vector bundle } \iff \Hom(R_\mu^\rk,E)=0 \]
for all coherent sheaves $E$ of rank $\rk$ and slope $\mu$.
\end{definition}\\

Raynaud's bundle $P_{2g+1}$ is a first example of a Raynaud bundle.
Theorem \ref{intro-1} could also be formulated as:
$P_{2g+1}$ is a Raynaud bundle $R^2_{g-1}$. 
We derive from this Theorem the existence Raynaud bundles $R^2_\mu$ for
all integer slopes $\mu$.
The main result of this paper is:

\begin{theorem}
For all pairs $(\rk,\mu)$ there exists a Raynaud bundle $R^\rk_\mu$.
\end{theorem}
This theorem is the equivalence (i)$\iff$(v) of Theorem \ref{t-1}. We
remark that such a Raynaud bundle is not unique. Indeed, twisting a
Raynaud bundle with
a line bundle of degree zero gives another, as well as taking the direct
sum of two such bundles. In Section \ref{CON} we construct the Raynaud
bundle $R^\rk_\mu$. Implicitly this construction appears in
Proposition 2 of \cite{HP}. However, there its construction is embedded
in the theory of the derived category. Here we work out this
construction, give the numerical invariants (Proposition \ref{prop-1}),
show the relation to base points of the $\Theta$-divisor, and give the
main properties in Theorem \ref{t-1}.

The purpose of section \ref{g-1-1} consists in  a fine tuning the
construction for the case when $\mu=g-1$. This allows the construction
of Raynaud bundles of smaller rank than those obtained in Section
\ref{CON}.  This way we obtain upper bounds
for $r$ for the base point freeness of $|R \cdot \Theta|$ on the moduli
spaces $\U_X(r,r(g-1))$ see Proposition \ref{prop-7} and Corollary
\ref{cor-4}. They imply upper bounds for the base point freeness
of the $\Theta$-divisor on $\SU_X(r)$ (see Proposition \ref{prop-8}).
However the bounds for base points of $|\Theta|$ on $\SU_X(r)$ are not
optimal see Arcara's result in \cite{Arc}
(see also older results of Popa in
\cite{Pop2}). O. Schneider used Raynaud's original bundles to produce
base points of $|\Theta|$ on $\SU_X(r)$ as extensions of Raynaud's bundle
by line bundles in \cite{Sch}.

For $X = \pdopi$ there exist semistable bundles $E$ only for integer
slopes. We see that the line bundle $\Ocal_\pdopi(\mu+1)$ is a Raynaud
bundle $R^\rk_\mu$ in this case. For an elliptic curve $X$ the existence
of Raynaud bundles is well known too (see Lemma 5 in \cite{HP1}). For
example: every
stable bundle $F$ of rank $\rk+1$ and degree one is a $R^\rk_0$.
Therefore, we may assume $g \geq 2$. 

\section{The Raynaud bundle $R_\mu^\rk$}\label{CON}
\subsection{Construction of $S_{\mu,R,m}$ for
$\mu \in [-g-1,-g)$}\label{cons-spez-mu}
Let $X$ be a smooth projective curve of genus $g$ over $k$.
We fix be a linebundle $L_1$ on $X$ of degree one.
By $L_{-1}$ we denote its dual.
Let $\mu=\frac{d}{r} \in [-g-1,-g)$ be a rational
number, where $d$ and $r$ are coprime integers with $r \geq 1$.
Furthermore, we fix a positive integer $R$.

A semistable vector bundle $E$ of slope $\mu$ is by definition a base
point of the linear system $|R \cdot \Theta |$, if for all vector
bundles $F$ with $\rk(F) =r \cdot R$ and $\det(F)=L_1^{\otimes
rR(g-1)-dR}$ we have $H^*(X, E \otimes F) \ne 0$.
(See also Beauville's article \cite{Bea} for the definition of base
points.)

%

We consider the two sheaves
\[ M_1=M_1(\mu,R):=
L_{-1}^{\otimes rR(g-1-\mu)}
\quad \mbox{and} \quad  
M_0=M_0(\mu,R):=\Ocal_X^{\oplus rR+1}.\]

We are interested in $M_1$ and $M_0$ because of the following lemma.

\begin{lemma}\label{cond-1}
Let $E$ be a semistable vector bundle of slope $\mu(E) \in [-g-1,-g)$.
Furthermore, we fix an integer $R \geq 1$, and the vector bundles $M_1$
and $M_0$ as above.
Then the following three conditions are equivalent:\\
\begin{tabular}{ll}
(i) & $E$ is no base point of $|R \cdot \Theta|$.\\
(ii) & For some morphism $M_1 \rarpa{\phi} M_0$
we have $H^0(E \otimes \coker(\phi)) = 0$.\\
(iii) & For some morphism $M_1 \rarpa{\phi} M_0$
the resulting morphism\\
& $H^1(E \otimes M_1) \to H^1(E \otimes M_0)$ is
injective.\\
\end{tabular}
\end{lemma}
\begin{proof}
Suppose there exists a vector bundle $F$ such that $H^*(E \otimes F)=0$
with $\rk(F)= rR$, and $\det(F) \cong L_1^{\otimes rR(g-1)-dR}$.
This implies that $F$ itself is a semistable bundle of slope
$\mu(F)=(g-1)-\mu$. Hence,
$F$ is semistable of slope
greater than $2g-1$. Thus, $F$
is globally generated. Indeed, we can generate this vector
bundle by $rR+1$ global sections. This way, we obtain a surjection:
$\Ocal^{\oplus rR+1} \rarpa{\pi} F$.
The kernel of $\pi$ is a line bundle $M_1$.
The determinant of this linebundle  is isomorphic to
$M_1 \cong \det(M_0) \otimes \det(F)\inv \cong \det(F)\inv$.
Thus, we have show, that (i) $\implies$ (ii).

To see the implication (ii) $\implies$ (i), we remark that setting
$F:=\coker(\phi)$ we obtain a sheaf of rank $rR$ and determinant
$L_1^{\otimes rR(g-1)-dR}$. From Riemann-Roch we deduce that $\chi(E
\otimes F)=0$. Thus, $H^0(E \otimes F)=0$ implies $H^*(E \otimes F)=0$.

The equivalence of (ii) and (iii) follows directly from the exact
sequence
\[H^0(E \otimes M_0) \to H^0(E \otimes \coker(\phi)) \to
H^1(E \otimes M_1) \to H^1(E \otimes M_0)\]
and the fact that $H^0(E \otimes M_0)=0$,
because the semistable bundle $E \otimes M_0$ has negative slope.
\qed
\end{proof}

We consider the vector space $V:=\Hom(M_1,M_0)$.
Its dimension is $v:=\dim(V)=(rR+1)(1-g+rR(g-1-\mu))$.
We consider the product space
\[\xymatrix{X&X \times \pdop(V\dual) \ar[l]_-p \ar[r]^-q
&\pdop(V\dual)}\]
Combining the universal morphism $M_1 \to V\dual \otimes M_0$ on $X$ and
the morphisms $\Ocal_{\pdop(V\dual)} \otimes V\dual \to
\Ocal_{\pdop(V\dual)}(1)$ on $\pdop(V\dual)$ we obtain a morphism
\[ \xymatrix{p^*M_1 \ar[r]^-\alpha
& p^*M_0 \otimes q^*\Ocal_{\pdop(V\dual)}(1)} \quad \mbox{on }X \times
\pdop(V\dual) .
\]
If we consider $\pdop(V\dual)$ as the moduli space of morphisms
different from zero from
$M_1(\mu,R)$ to $M_0(\mu,R)$ modulo scalars, then
this morphism is the universal family over $\pdop(V\dual)$.
Since $p^*M_1(\mu,R)$ is a linebundle and $\alpha$ is not trivial, we deduce
that $\alpha $ is injective. We denote its cokernel by $G(\mu,R)$.
Twisting the the short exact sequence
\[ 0 \to p^*M_1(\mu,R) \to p^*M_0(\mu,R) \otimes q^*\Ocal_{\pdop(V\dual)}(1)
\to G(\mu,R) \to 0\]
by $q^*\Ocal_{\pdop(V\dual)}(m)$ for any $m \geq0$,
we obtain a short exact
sequence of sheaves possessing no higher direct images with respect to
$p$.
Thus, we obtain a short exact sequence on $X$.
\[ 0 \to \Sym^m(V\dual)\otimes M_1(\mu,R) \to
 \Sym^{m+1}(V\dual)\otimes M_0(\mu,R) \to
p_*(G(\mu,R) \otimes q^*\Ocal_{\pdop(V\dual)}(m)) \to 0 \]

We define the sheaf $S_{\mu,R,m}:=p_*(G(\mu,R) \otimes
q^*\Ocal_{\pdop(V\dual)}(m))$.
From the construction of $S_{\mu,R,m}$ we conclude the following
the first properties of the sheaf $S_{\mu,R,m}$, namely

\begin{proposition}\label{prop-0}
There exists a short exact sequence
\[ 0 \to \Sym^m(V\dual)\otimes M_1(\mu,R) \to
\Sym^{m+1}(V\dual)\otimes M_0(\mu,R) \to S_{\mu,R,m} \to 0 .\]
\end{proposition}

\begin{proposition}\label{prop-1}
The numerical invariants of the sheaf $S_{\mu,R,m}$ are given by\\
\[\begin{array}{rcl}
\rk(S_{\mu,R,m})& = &\binom{m+v-1}{m}\left( (rR+1)\frac{v+m}{m+1}-1 \right)\\
\deg(S_{\mu,R,m}) & =&
\binom{m+v-1}{m}rR(g-1-\mu)\\
\mu(S_{\mu,R,m}) & = &\frac{(m+1)rR(g-1-\mu)}{(m+1)rR+(v-1)(rR+1)} =
g-1-\mu -\frac{(v-1)(rR+1)(g-1-\mu)}{(m+1)rR+(v-1)(rR+1)}        \\
\end{array}\]
\end{proposition}

\begin{remark}
Considered as a function depending on $m \in \ndop$ the slope of
$S_{\mu,R,m}$ is of the form
$\mu(S_{\mu,R,m})=g-1-\mu-\frac{a}{m+b}$ for some positive $a,b \in
\qdop$.
\end{remark}

\subsection{Properties of $S_{\mu,R,m}$ for
$\mu \in [-g-1,-g)$}\label{prop-spez-mu}

We keep the notation of \ref{cons-spez-mu}. In particular we still
assume that $\mu =\frac{d}{r} \in [-g-1,-g)$, $R$ and the bundles $M_1$
and $M_0$ are fixed in \S \ref{prop-spez-mu}. 
We need the following result.

\begin{lemma}\label{lemma-1}(\cite[Lemma 13]{HP})
Let $U$ and $W$ be $k$-vector spaces of finite dimension.
Suppose that a given morphism $\xymatrix{ U \otimes \Ocal_{\pdop^n}
\ar[r]^-\rho & W \otimes \Ocal_{\pdop^n}(1)}$ on $\pdop^n$
is not injective.
Then for any integer
$m \geq (\dim(U)-1)n$ we have $H^0(\ker(\rho)(m)) \ne 0$.
\end{lemma}

Note, that the sheaf $E$ in the following proposition is not necessarily
of slope $\mu$. However, semistable vector bundles of negative slope
fulfill the premise of the proposition.

\begin{proposition}\label{prop-2}
Let $E$ be a sheaf on $X$ with the property that $H^0(X,E)=0$.
For any $m \geq (v-1)(h^1(M_1 \otimes E)-1)$,
the following two conditions are equivalent:\\
\begin{tabular}{ll}
(i) & There exists a short exact sequence
$0 \to M_1 \to M_0 \to F \to 0$ with $H^0(F \otimes E )=0$.\\
(ii) & $H^0(S_{\mu,R,m} \otimes E) = 0$.\\
\end{tabular}
\end{proposition}
\begin{proof}
From $H^0(E)=0$, we deduce that $E$ is a vector bundle
and the equalities $h^0(M_0 \otimes E)=0=h^0(M_1 \otimes E)$.
Furthermore, the dimension
$h^1(M_1 \otimes E)$ can be computed by Riemann-Roch and is positive.

We consider the short exact sequence
\[ 0 \to p^*(M_1 \otimes E) \to p^*(M_0 \otimes E) \otimes
q^*\Ocal_{\pdop(V\dual)}(1) \to p^*E \otimes G(\mu,R) \to 0\]
on $X \times \pdop(V\dual)$. 
Since $H^0(M_0 \otimes E)=0$ we obtain on $\pdop(V\dual)$ an exact
sequence
\[0 \to q_*(p^*E \otimes G(\mu,R) \to H^1(M_1 \otimes E) \otimes
\Ocal_{\pdop(V\dual)} \rarpa{\beta} H^1(M_0 \otimes E) \otimes
\Ocal_{\pdop(V\dual)}(1) \to\]
By base change, condition (i) is equivalent to the injectivity of the
morphism $\beta$. This is by lemma \ref{lemma-1} equivalent to
$H^0(q_*(p^*E \otimes G(\mu,R)) \otimes \Ocal_{\pdop(V\dual)}(m))=0$.
Thus, (i) is equivalent to $H^0(p^*E \otimes G(\mu,R) \otimes 
q^*\Ocal_{\pdop(V\dual)}(m))=0$.
This implies the result by definition of $S_{\mu,R,m}$, and the
projection formula.
\qed
\end{proof}

As a corollary of the proof of proposition \ref{prop-2} we obtain
the
\begin{corollary}\label{cor-1}
For any sheaf $E$ on $X$ the assignment
$m \mapsto h^0(S_{\mu,R,m} \otimes E)$ is the Hilbert function of the
torsion free sheaf $q_*(p^*E \otimes G(\mu,R))$.
In particular, $h^0(S_{\mu,R,m} \otimes E) \ne 0$ implies
$h^0(S_{\mu,R,M} \otimes E) \ne 0$ for all $M \geq m$.
\end{corollary}

\begin{corollary}\label{cor-2}
For any $m \geq 0$ the sheaf $S_{\mu,R,m}$ is a vector bundle on $X$.
\end{corollary}
\begin{proof}
Take a stable vector bundle $F$ with
$\det(F) \cong L_1^{\otimes rR(g-1-\mu)}$. As seen in Lemma
\ref{cond-1}, there exists a short exact sequence
$0 \to M_1 \rarpa{\phi} M_0 \rarpa{\pi} F \to 0$.
We take a (sufficiently negative) linebundle $E$ on $X$,
such that $h^0(E)=0=h^0(E \otimes F)$.
The line bundle $E$ fulfills the assumption of
Proposition \ref{prop-2} and condition (i) of Proposition \ref{prop-2}
is satisfied. Thus, we conclude $H^0(S_{\mu,R,M} \otimes E) =0$ for $M
\gg 0$. By Corollary \ref{cor-1} this implies
$H^0(S_{\mu,R,m} \otimes E) =0$ for all $m \geq 0$.
Hence, $S_{\mu,R,m}$ is torsion free.
\qed
\end{proof}

\subsection{Definition and properties of $S_{\mu,R}^{\rk}$ for
$\mu \in [-g-1,-g)$}\label{prop-spez-mu-r}
In this part \ref{prop-spez-mu-r} we still assume that $\mu=\frac{d}{r}
\in [-g-1,-g)$ and $d$ and $r$ are coprime. Thus, for a vector bundle $E$
of slope $\mu(E)=\mu$ we have $\rk(E):=h \cdot r$ for some natural
number $h$. Remember, the number $v=(rR+1)(1-g+rR(g-1-\mu))$.
For any number $\rk$ which is  a multiple of $r$ we define
\[S_{\mu,R}^{\rk} := S_{\mu,R,(v-1)(\rk \cdot (g-1-\mu)(rR+1)-1)} \,.\]
\begin{proposition}\label{prop-3}
For a semistable vector bundle $E$ of slope $\mu(E) = \mu \in [-g-1,-g)$
on the curve $X$ we have an equivalence.
\[ E \mbox{ is not a base point of } |R \cdot \Theta| \iff
H^0(S_{\mu,R}^{\rk(E)} \otimes E) = 0 .\]
\end{proposition}
\begin{proof}
First we note that $h^0(E)=0$ because $E$ is semistable of negative
slope. Thus, $h^0(M_1 \otimes E)=0$ and we can compute
$h^1(M_1 \otimes E)=\rk(E) \cdot (g-1-\mu)(rR+1)$ by the Riemann-Roch
theorem. Now we deduce the equivalence from Propositions \ref{cond-1}
and \ref{prop-2} because we took the number $m$ in the definition of
$S_{\mu,R}^{\rk}$ to be the smallest possible $m$ such that
Proposition \ref{prop-2} applies to $S_{\mu,R,m}$ and $E$. \qed
\end{proof}

\begin{lemma}\label{lemma-2}
If $E$ is a coherent sheaf of slope $\mu(E) = \mu \in [-g-1,-g)$ on $X$
with the property $H^0(S_{\mu,R}^{\rk(E)} \otimes E) = 0$,
then $E$ is semistable.
\end{lemma}
\begin{proof}
First, we note that $H^0(S_{\mu,R}^{\rk(E)} \otimes E) = 0$ implies that
$E$ is torsion free. Now suppose that $E$ is not semistable. Let $E'
\subset E$ be a destabilizing subbundle. We have $\mu(E') \geq
\mu(E)+\frac{1}{\rk(E)(\rk(E)-1)}$.
By proposition \ref{prop-1} and the choice of $m$ we derive the
inequality
\[ \mu(S_{\mu,R}^{\rk(E)} \otimes E') = \mu(S_{\mu,R}^{\rk(E)}) +
\mu(E') > g-1 .\]
This implies $\chi(S_{\mu,R}^{\rk(E)} \otimes E') >0$. Hence, we deduce
$0 \ne H^0(S_{\mu,R}^{\rk(E)} \otimes E') \subset H^0(S_{\mu,R}^{\rk(E)}
\otimes E) $, which contradicts our assumption.
\qed
\end{proof}

\subsection{Definition and properties of $S_\mu^{\rk}$ for $\mu \in
[-g-1,-g)$}\label{prop-spez-rk-mu}
We define the vector bundle $S_\mu^\rk$ to be $S_{\mu,\tilde R}^\rk$
with $\tilde R =\lceil \frac{(\rk+1)^2}{4} \rceil$.
Still assuming, that $\mu \in [-g-1,-g)$, $\mu=\frac{d}{r}$, with $\rk=r
h$ for some integer $h$ we conclude the following result.

\begin{proposition}\label{prop-5}
For a coherent sheaf $E$ of slope $\mu \in [-g-1,-g)$
and rank $\rk$, we have the equivalence
\[E \mbox{ is a semistable vector bundle} \iff
H^0( S_\mu^\rk \otimes E)=0 .\]\end{proposition} 
\begin{proof}
We have seen in Lemma \ref{lemma-2}, $h^0(E \otimes S_\mu^\rk)=0$
implies that $E$ is a semistable vector bundle. Suppose now that $E$ is
a semistable vector bundle. Since the generalized $\Theta$-divisor $|R
\cdot \Theta|$ is base point free for $R \geq \frac{(\rk+1)^2}{4}$
(see Theorem 4.1 in
Popa's article \cite{Pop}) we have by Proposition \ref{prop-3}, that
$h^0(S^{\rk}_{\mu,R} \otimes E)=0$ for all $R \geq \frac{(\rk+1)^2}{4}$.
By definition of $S_\mu^\rk$ this proves the claimed statement.
\qed
\end{proof}

\subsection{The vector bundles $S_\mu^\rk$ and $R_\mu^\rk$}
Let $\mu=\frac{d}{r}$ be a rational number expressed as quotient of
two coprime integers with $r \geq 1$.
In contrast to parts \ref{cons-spez-mu}--\ref{prop-spez-rk-mu} there
exists no restriction on $\mu$. We take an integer $\rk$ which is a
multiple of $r$. 

We define the vector bundle $S_\mu^\rk$ through
\[S_\mu^\rk :=
L_{-1}^{\otimes (\lfloor \mu \rfloor + 1+g)} \otimes
S_{\mu -(\lfloor \mu \rfloor + 1+g)}^\rk \, .\]
This is well defines, as $\mu -(\lfloor \mu \rfloor + 1+g) \in
[-g-1,-g)$. 
Now, we have the
\begin{proposition}\label{prop-6}
If $E$ is a coherent sheaf of positive rank $\rk$ and of slope $\mu$,
then
\[E \mbox{ is a semistable vector bundle} \iff
H^0( S_\mu^\rk \otimes E)=0 .\]\end{proposition}
\begin{proof}
We have $E$ is a semistable vector bundle,
if and only if $E \otimes L_{-1}^{\otimes
(\lfloor \mu \rfloor + 1+g)}$ is a semistable vector bundle. 
Since $\mu( L_{-1}^{\otimes
(\lfloor \mu \rfloor + 1+g)} \otimes E) = \mu-\lfloor \mu \rfloor -g-1 \in
[-g-1,-g)$, we can apply Proposition \ref{prop-5} to obtain that
$E$ is a semistable vector bundle, if and only if the cohomology group
$H^0(S^\rk_{\mu-(\lfloor \mu \rfloor +g+1)}
\otimes L_{-1}^{\otimes
(\lfloor \mu \rfloor + 1+g)} \otimes E)$ is zero. By definition this
happens exactly when $H^0( S_\mu^\rk \otimes E)=0$.
\qed
\end{proof}

We define the vector bundle $R^\rk_\mu$ to be the dual of $S_{\mu}^\rk$.
We have the

\begin{theorem}\label{t-1}
Let $E$ be a coherent sheaf on the smooth projective curve $X$ of rank
$\rk >0$ and slope $\mu=\frac{d}{r}$.
The following conditions are equivalent:\\
\begin{tabular}{ll}
(i) & E is a semistable vector bundle.\\
(ii) & There exists a sheaf $F$ of rank $\lceil  \frac{(\rk+1)^2}{4}
\rceil r$ such that $H^*(E \otimes F)=0$.\\
(iii) & There exists a sheaf $F\ne 0$ such that $H^*(E \otimes
F)=0$.\\
(iv) & $H^0(S^\rk_\mu \otimes E)=0$.\\
(v) & $\Hom(R_{\mu}^\rk,E) = 0$.\\
\end{tabular}
\end{theorem}
\begin{proof}
The implications (ii) $\implies$ (iii) $\implies$ (i), and (iv) $\iff$ (v)
are standard.
The equivalence of (i) and (ii) is shown in Theorem 4.1 of \cite{Pop}.
The equivalence of (i) and (iv) was shown in Proposition \ref{prop-6}.
\qed
\end{proof}

\subsection{Further remarks}
Let $R_\mu^\rk$ be a Raynaud bundle constructed above.
We remark that for any unstable $E$ of slope $\mu$ and rank $\rk$ we
have $\hom(R_\mu^\rk,E')-\ext^1(R_\mu^\rk,E') > 0$ for all destabilizing
$E' \subset E$ (see Lemma \ref{lemma-2}). Suppose $R_\mu^\rk$ is not
stable. Then we have a surjection to a stable bundle
$R_\mu^\rk \to \ol{R_\mu^\rk}$ with
$\mu(R_\mu^\rk) \geq \mu( \ol{R_\mu^\rk})$. The last inequality implies
$\hom(\ol{R_\mu^\rk},E')-\ext^1(\ol{R_\mu^\rk},E') > 0$ for all $E'$ as above.
Since $\Hom(R_\mu^\rk,E)=0$ for all semistable $E$, we deduce that
$\Hom(\ol{R_\mu^\rk},E)=0$. As a consequence we note:

\begin{proposition} 
There are stable Raynaud bundles $R_\mu^\rk$.
\end{proposition}

\begin{remark} The semicontinuity Theorem (III.12.8 in \cite{Har})
implies that semistability is an open condition. Indeed, take any vector
bundle $R$ and define $R$-semistability of $E$ by the condition
$\Hom(R,E)=0$. From the semicontinuity Theorem we deduce that
$R$-semistability is an open condition in flat families.
\end{remark}

\begin{question}
What is the smallest possible rank for a Raynaud bundle $R^\rk_\mu$?
As we see in Section \ref{g-1-1}, there can be constructed Raynaud
bundles of smaller rank. However, these bundles still have huge rank as
we can see in the small table after Corollary \ref{cor-3}. It is the
author's  believe that these ranks are still far from being optimal.
\end{question}

\section{Base points of $|R \cdot \Theta|$ on
$\U_X(r,r(g-1))$}\label{g-1-1}
Throughout this section \ref{g-1-1} $E$ is a coherent sheaf of rank $r$ and slope
$\mu(E)=g-1$. It turns out that in this case we can construct vector
bundles $S^r_R(M_0)$ with the same property like the bundle $S^r_{(g-1),R}$
given in Proposition \ref{prop-3} having a significant smaller rank
than $S^r_{(g-1),R}$.

\subsection{A Raynaud bundle for $\mu=g-1$}\label{g-1-1-1}
Let us fix the notation: We consider a smooth projective curve $X$ of
genus $g \geq 2$ over an algebraically closed field $k$.
Furthermore, we fix a natural number $R \geq 2$.

\begin{lemma}\label{nice-triple}
There exists a short exact sequence of vector bundles on $X$
\[0 \to M_1 \rarpa{\phi} M_0 \rarpa{\psi} F \to 0\]
with the  following properties:\\
\begin{tabular}{lp{14cm}}
(i) & $F$ is stable with $\rk(F)=R$, and $\det(F) \cong \Ocal_X$.\\
(ii) & $M_0$ is stable with $\rk(M_0)=R+1$,
and $\deg(M_0)=(R+1)(1-g)-1$.\\
(iii) & $\Ext^1(M_0,F)=0$.\\
\end{tabular}
\end{lemma}
\begin{proof}
Considering all triples $M_0 \rarpa{\psi} F$ we see that there exist
surjections $M_0 \rarpa{\psi} F$ with the given numerical invariants and
$F$ stable (see Proposition 7.3 and Theorem 7.7 in \cite{BGG}).
Take a pair $(\tilde M_0,F)$ of stable sheaves with $\det(F) \cong
\Ocal_X$. $\rk(F)=R$, $\deg(\tilde M_0)=(R+1)(1-g)$,
$\rk(\tilde M_0)=R+1$ such that $H^*(F \otimes \tilde M_0\dual)=0$.
The existence of such a pair is well known (cf. Beauville's survey
article \cite{Bea}). 
The stability of $\tilde M_0$, and $\mu(\tilde M_0) \in \zdop$ imply
that for any surjection $\pi: \tilde M_0 \to k(P)$ the kernel $M_0$ is
also stable. From the short exact sequence 
$0 \to M_0 \to \tilde M_0 \to k(P) \to 0$, and $H^*(F \otimes \tilde
M_0\dual)=0$ we deduce that $M_0$ satisfies (ii) and (iii).
Since the properties (i)--(iii) are open properties on the irreducible
moduli space of triples $M_0 \rarpa{\psi} F$ (again Theorem 7.7 in
\cite{BGG}) we deduce the claim.
\qed \end{proof}

\begin{notation}
From now on we take fixed vector bundles $M_1$ and $M_0$ from a short
exact sequence $0 \to M_1 \rarpa{\phi} M_0 \rarpa{\psi} F \to 0$ like in
Lemma \ref{nice-triple}.
Compare the following result with Lemma \ref{cond-1}.

\end{notation}
\begin{lemma}\label{cond-1-ft}
Let $E$ be a semistable vector bundle of slope $\mu(E) = g-1$.
Furthermore, we fix an integer $R \geq 2$, and the vector bundles $M_1$
and $M_0$ as above.
Then the following three conditions are equivalent:\\
\begin{tabular}{ll}
(i) & $E$ is no base point of $|R \cdot \Theta|$.\\
(ii) & For some morphism $M_1 \rarpa{\phi} M_0$
we have $H^0(E \otimes \coker(\phi)) = 0$.\\
(iii) & For some morphism $M_1 \rarpa{\phi} M_0$
the resulting morphism\\
& $H^1(E \otimes M_1) \to H^1(E \otimes M_0)$ is
injective.\\
\end{tabular}
\end{lemma}
\begin{proof}
The implications (iii) $\iff$ (ii) $\implies$ (i) follow like in the
proof of \ref{cond-1}. The problem with (i) $\implies$ (ii) is that not
all semistable vector bundles $F$ of rank $R$ and determinant $\Ocal_X$
are quotients of $M_0$.
Applying $\Hom( - ,F)$ to the short exact sequence of lemma
\ref{nice-triple} yields the long exact sequence
\[ \to \Hom( M_1,F) \rarpa{\alpha} \Ext^1(F,F) \to
\Ext^1(M_0,F) \to \Ext^1(M_1,F) \to 0\]
The consequences of the vanishing of $\Ext^1(M_0,F)$ (see Lemma
\ref{nice-triple} (iii))
we express in terms of the Quot scheme $\Quot = \Quot(M_0)_X^{R,0}$  
of rank $R$, degree zero quotients of $M_0$.\\
First we conclude, that $\Ext^1(M_1,F)=0$. This is the
obstruction space of the Quot scheme at $[\psi]=[\psi:M_0 \to F]
\in \Quot(k)$. Thus, there exists a smooth open
neighborhood $U$ of $[\pi]$ which parameterizes semistable vector
bundles.\\
Secondly we deduce surjectivity of $\alpha$. This is the tangent map at
$[\psi]$ of the morphism $U \to \U_X(R,0)$ from $U$ to the moduli space
of rank $R$ bundles of degree zero. Passing to a smaller open subset of
$U$ we may assume that $U \to \U_X(R,0)$ is a smooth morphism. The
image $V$ of $U$ is open and contains a vector bundle with trivial
determinant. We conclude, that a dense open subset $V_{\Ocal_X}$ of the
moduli space $\SU_X(R,\Ocal_X)$ of rank $R$ bundles with trivial
determinant is parameterized by points of our Quot scheme.\\
Now assume (i). Thus, there exists a vector bundle $F$ of rank $R$ with
trivial determinant, such that $H^*(X,E \otimes F)=0$.
Thus, the vector bundles $G$ parameterized by $\SU_X(R,\Ocal_X)$ with
$h^1(E \otimes G) \ne 0$ form a divisor which can not contain the open
set $V_{\Ocal_X}$. This shows that (ii) holds.
\qed
\end{proof}

Now we set $V:=\Hom(M_1,M_0)$. Since the difference of the slopes
$\mu(M_0)-\mu(M_1) > 2g-2$, we have $\Ext^1(M_1,M_0)=0$ and can compute
the dimension $v$ of $\Hom(M_1,M_0)$ by Riemann-Roch to be
$v=(R+1)(R-1)(g-1)+R$.\\ 
We follow the construction in 2.1:
We consider the projections $\xymatrix{X&X \times \pdop(V\dual)
\ar[l]_-p \ar[r]^-q &\pdop(V\dual)}$ and the morphism
$\alpha: p^*M_1 \to p^*M_0 \otimes \Ocal_{\pdop(V\dual)}(1)$ 
to obtain for any $m \geq 0$ the bundle $S_{R,m}(M_0):=p_*( \coker(\alpha)
\otimes q^*\Ocal_{\pdop(V\dual)}(m))$.
We set
\[S_R^r(M_0):=S_{R,w} \quad \mbox{ with }
w:=((R-1)(R+1)(g-1)+R-1)(r(R+1)(g-1)+r-1) \, .\]
\[\mbox{and } S^r(M_0):=S^r_u(M_0) \quad \mbox{ with }
u:=\left\lceil \frac{(r+1)^2}{4} \right\rceil 
\mbox{, and } R^r(M_0)=\left( S^r(M_0) \right)\dual.\]

\begin{theorem}\label{prop-g-1}
{\bf (Properties of $S_{R,m}(M_0)$, $S_R^r(M_0)$, and $S^r(M_0)$) }\\
\begin{tabular}{lp{14.5cm}}
(i) & $S_{R,m}(M_0)$, $S_R^r(M_0)$, and $S^r(M_0)$ are vector bundles on
$X$.\\
(ii) & The numerical invariants of $S_{R,m}(M_0)$ are given by
\[ \begin{array}{rcl}
\deg(S_{R,m}(M_0)) & = &  ((R+1)(1-g)-1)\frac{v-1}{m+1}\binom{v+m-1}{m}\\
\rk(S_{R,m}(M_0)) & = &  \frac{Rv+Rm+v-1}{m+1} \binom{v+m-1}{m}\\
\mu(S_{R,m}(M_0)) & = & \frac{((R+1)(1-g)-1)(v-1)}{Rv+Rm+v-1}
\end{array} \] 
where $v:=(R+1)(R-1)(g-1)+R$.\\
(iii) & For $m \geq 0$, and any coherent sheaf $E$ on $X$ we have
$H^0(S_{R,m}(M_0) \otimes E) \ne 0$ implies
$H^0(S_{R,M}(M_0) \otimes E) \ne 0$ for all $M \geq m$.\\
(iv)&For a semistable sheaf $E$ of rank $r$ with $\chi(E)=0$ we have an
equivalence
\[ E \mbox{ is a base point of } |R \cdot \Theta| \mbox{ on } \U_X(r,r(g-1))
\iff H^0(S_R^r(M_0) \otimes E) \ne 0 .\]\\
(v)&For a coherent sheaf $E$ of rank $r$ with $\chi(E)=0$ we have an
equivalence
\[ E \mbox{ is semistable }
\iff H^0(S^r(M_0) \otimes E)=0 \iff \Hom(R^r(M_0),E) =  0 .\]\\
\end{tabular}
\end{theorem}
\begin{proof}
The results follow straightforward by applying Lemma \ref{cond-1-ft}
instead of Lemma \ref{cond-1}. In particular:
(i) follows from Corollary \ref{cor-2}, (ii) from Proposition
\ref{prop-1}, (iii) is Corollary \ref{cor-1}, Proposition \ref{prop-3}
gives (iv), and (v) is just Proposition \ref{prop-6}. 
\qed
\end{proof}

\begin{corollary}\label{cor-3}
The slope $\mu(S_{R,m}(M_0))$ of the vector bundle $S_{R,m}(M_0)$
considered as a function of $m$ is
of type $\mu(S_{R,m}(M_0)) =
\frac{-a}{Rm+b}$ for positive integers $a,b
\in \ndop$. In particular, we have
\[ \mu(S_{R,m}(M_0)) \geq 1-g
\iff m \geq (R-1)+\frac{R-g}{R(g-1)} \,.\]
\end{corollary}

We list the rank and the slopes of the Raynaud bundles
$R^r(M_0)$ which we obtained for $\mu=g-1$ by the methods of this
subsection for $r,g \in \{2,3,4\}$.\\
\begin{tabular}{r|r|r|r}
$g$&$r$& $\rk(R^r(M_0))$ & $\mu(R^r(M_0))$ \\
\hline \hline
2 & 2 &59539855602920 & $\frac{50}{313}$\\
\hline
2 & 3 &641752198359834620231606142864&$\frac{54}{659}$\\
\hline
2 & 4 & $5.78978673052\cdot 10^{106}$ & $\frac{486}{13669}$\\
\hline
3 & 2 & 483505260221028663042477162264 & $\frac{54}{331}$\\
\hline
3 & 3 & $4.88907844550 \cdot 10^{63}$ & $\frac{363}{4393}$\\
\hline
3 & 4 & $2.18037666849 \cdot 10^{230}$ & $\frac{1734}{48661}$\\
\hline
4 & 2 & 182463883365641199732269260672875437828878976664 &
$\frac{338}{2057}$\\
\hline
4 & 3 & $5.06529456824 \cdot 10^{100}$ & $\frac{192}{2317}$\\
\hline
4 & 4 & $1.52141697065 \cdot 10^{364}$ & $\frac{3750}{105157}$\\
\end{tabular}

These values show that even for small values of $g$ and $r$ the help a
computer program ({\tt bc} in my case) is needed to compute the rank and
slopes of the Raynaud bundles.

\subsection{Base points of $|2 \cdot \Theta|$ on $\U_X(r,r(g-1))$ and of
$|\Theta|$ on $\SU_x(r)$}
\begin{lemma}\label{lem-1}
Let $F$ be a vector bundle of rank $r_F$ and slope $\mu(F) \leq g-1$. If
$r_E \geq r_F$ is an integer, then there
exists a semistable vector bundle $E$ of rank $r_E$ and slope
$\mu(E)=g-1$ with $\Hom(F,E) \ne 0$. Moreover, if $r_E > r_F$ or
$\mu(F)<g-1$, then the S-equivalence classes of the bundles $E$ with
$\Hom(F,E) \ne 0$ form a
positive dimensional subset in the moduli space $\U_X(r_E,r_E(g-1))$.
\end{lemma}
\begin{proof}
The proof works by induction on the rank $r_F$. We take an elementary
transformation $0 \to F \to \tilde F \to T \to 0$ such that $\tilde F$
is a vector bundle of rank $r_F$ and $\mu(\tilde F)=g-1$. Now we
distinguish two cases:

{\bf Case 1:} $\tilde F$ is stable. In this case we may take $E=\tilde F
\oplus E'$ to be a sum of two stable vector bundles of slope $g-1$.

{\bf Case 2:} If $\tilde F$ is not stable there exists a surjection
$\tilde F \to F''$ to a stable bundle $F''$ with $\mu(F'') \leq g-1$ and
$\rk(F'')<\rk(F)$. Thus, by the induction hypothesis we are done.

We remark that for $r_F=1$ we are always in the situation of case 1.
The statement about the dimensions is trivial (we may change the
determinant of $\tilde F$
by varying the support of $T$ or vary the bundle $E'$).  \qed
\end{proof}

\begin{proposition}\label{prop-7}
For any smooth projective curve $X$ of genus $g \geq 2$ the linear
system $|2 \cdot \Theta|$ on the moduli space $\U_X(r,r(g-1))$ has base
points for $r \geq \frac{27g^2-15g+2}{2}$. Furthermore, the base locus
is of positive dimension for $g>2$ or $r > \frac{27g^2-15g+2}{2}$
\end{proposition}
\begin{proof}
The dual vector bundle $(S_2,1(M_0))\dual$ has slope
$\mu((S_2,1(M_0))\dual) \leq g-1$ by Corollary \ref{cor-3} and is of
rank $\frac{27g^2-15g+2}{2}$ by Proposition \ref{prop-g-1}.(ii).
Thus, for all $r \geq \frac{27g^2-15g+2}{2}$
we find semistable vector bundles $E$ with $\Hom((S_2,1(M_0))\dual,E)
\ne 0$. This is equivalent to $H^0(S_2,1(M_0) \otimes E) \ne 0$ and
implies by (iii) and (iv) of Proposition \ref{prop-g-1} that $E$ is a
base point of $|2 \cdot \Theta|$.
\qed
\end{proof}

\begin{proposition}\label{prop-8}
For any smooth projective curve $X$ of genus $g \geq 2$ the linear
system $|\Theta|$ on $\SU_X(r)$ has base points for $r \geq
\frac{27g^2-15g+2}{2}$. The base locus is of positive dimension for
$r > \frac{27g^2-15g+2}{2}$.
\end{proposition}
\begin{proof}
We take a base point $[E] \in \U_X(r,r(g-1))$. There exists a line
bundle $M \in \Pic^{1-g}(X)$ such that $\det(E \otimes M) \cong
\Ocal_X$. We claim that $E \otimes M$ is a base point of $|\Theta |$ on
$\SU_X(r)$. Indeed if it were not a base point, we would have a proper
divisor $D \subset \Pic^{g-1}$ such that for all $L \in Pic^{g-1}(X)
\setminus D$ we have $H^*(E \otimes M \otimes L)=0$. Take $L \in
\Pic^{g-1}(X)$, such that neither $L$ nor $(M^{-2} \otimes L\inv)$ are in
$D$. Then it follows that $H^*(E \otimes M \otimes (L \oplus (M^{-2}
\otimes L\inv)))=0$. However, $\det(M \otimes (L \oplus (M^{-2}
\otimes L\inv))) \cong \Ocal_X$. Thus, would not be a base point. This
proves the claim.
\qed
\end{proof}

\subsection{Base points of $|R \cdot \Theta|$ on
$\U_X(r,r(g-1))$}\label{g-1-1-3}
As in the subsection before we remark that $S_{R,R}(M_0)$ has slope at
least $1-g$ (see Corollary \ref{cor-3}). Thus, by Lemma \ref{lem-1}
we obtain that for all $r \geq \rk(S_{R,R}(M_0))$ the linear system $|R
\cdot \Theta|$ is not base point free on $\U_X(r,r(g-1))$. Moreover, if
$g \geq R$, then this holds for all $r \geq \rk(S_{R,R-1}(M_0))$.\\
\begin{minipage}{7cm}
In the table we have computed for small $R$ and $g$ the minimal ranks
$r$ for which $\U_X(r,r(g-1))$ is known to have base points by our
method. 
\end{minipage}
\hspace{1em}
\begin{tabular}{r||r|r|r|r}
& g=2 & g=3 & g=4 & g=5\\
\hline
\hline
R=2 & 40 & 100 & 187 & 301\\
R=3 & 3718 & 5130 & 14238 & 30450\\
R=4 & 160930 & 2443665 & 1332800 & 3786640\\
\end{tabular}\\
The big values of $r$ explain why we do not include an explicit formula
in the next corollary. However, the interested reader can extract the
rank using Theorem \ref{prop-g-1} (ii).

\begin{corollary}\label{cor-4}
For any $R \geq 2$ there exists a polynomial $p_R$ of degree $R$ such
that for all $r \geq p_R(g)$ the linear system $|R \cdot \Theta|$ on
$\U_X(r,r(g-1))$ is not base point free.
\qed
\end{corollary}

\vfill
Georg Hein, Universit\"at Duisburg-Essen, Fachbereich Mathematik, 45117
Essen\\email: {\tt georg.hein@uni-due.de}
\end{document}